\documentclass[letterpaper,12pt,leqno,oneside]{article}
\usepackage{color}
\usepackage{bm}
\usepackage{exscale}
\usepackage{amsmath}
\usepackage{amsfonts}
\usepackage{stmaryrd}
\usepackage{amscd}
\usepackage{graphicx}
\usepackage{amsxtra}
\usepackage{amssymb}
\usepackage{theorem}
\usepackage[final]{epsfig}
\usepackage{eqnarray}
\usepackage{hyperref}
\newtheorem{proposition}{Proposition}[subsection]

{\theorembodyfont{\rmfamily}\newtheorem{remark}[proposition]{Remark}}

{\theorembodyfont{\rmfamily}}

\newfont{\abc}{cmtt10 scaled 1200}

\def\R{\mathbb{R}}

\def\Z{\mathbb{Z}}
\def\Q{\mathbb{Q}}

\def\ve{\varepsilon}

\def\ra{\rightarrow}
\def\cs{\symbol{35}}
\def\p{\partial}
\def\qed{\hfill $\Box$ \\}
\def\mm{\mbox}
\def\v{= \emptyset}
\def\n{\neq \emptyset}

\def\si{$\mathcal{S}$}

\def\bp{\langle A \rangle}

\begin{document}

\vspace*{0cm}

\begin{center}\Large{\bf{Contracting Maps and Scalar Curvature}}\\
\bigskip
\large{\bf{Joachim Lohkamp}\\
\medskip}

\end{center}

\noindent Mathematisches Institut, Universit\"at M\"unster, Einsteinstrasse 62, Germany\\
 {\small{\emph{e-mail: j.lohkamp@uni-muenster.de}}}

\setcounter{section}{1}
\renewcommand{\thesubsection}{\thesection}
\subsection{Introduction} \label{int}

In these short notes we explain how to derive the following size estimate in scalar curvature geometry for \emph{contracting maps}. A $C^1$-map $f:V\ra W$ between Riemannian manifolds $V$ and $W$  is \emph{$\ve$-contracting} when $\|Df(v)\| \leq \ve \cdot \|v\|$ for any tangent vector $v \in TV$.  $1$-contracting maps are called \emph{contracting maps}.\\

\textbf{Theorem} \textbf{(Contractions and Scalar Curvature)} {\itshape For $n \ge 2$, let $M^n$  be a compact manifold with boundary  $\partial M^n$ and  $f:(M^n, \partial M^n)\rightarrow(S^{n},\{p\})$ a $C^1$-map of non-zero degree mapping a neighborhood of $\p M^n$ to $p \in S^{n}$. Then, there is some $\sigma_n \ge 1$, that depends only on $n$, so that for any metric $g$ on $M^n$ making $f$ a contracting map:
\begin{equation}\label{csc}
scal(g)(z) \le \sigma_n, \mm{ in some point } z \in M^n,
\end{equation}
where $scal(g)$ denotes the scalar curvature of the smooth Riemannian metric $g$.}\\

This result was previously known in, what will henceforth be called,  the \emph{classical cases} where $M^n$ either admits a spin structure or it has dimension $n\le 7$  [GL1],Ch.12.
The dimensional constraint arises from the use of minimal hypersurfaces in the non-spin case. These hypersurfaces can carry difficult singularities, occurring in dimensions $\ge 8$, which could not be handled by classical means.\\

Our extension to the general case is a simple application of scalar curvature splitting theorems. (For the purposes of this paper the basic version derived in [L1], Th.1 and the naturality theory [L2], Th.3 are sufficient.) Roughly speaking, such splitting results shift  arguments involving (potentially singular) minimal hypersurfaces back into an entirely smooth scenario.\\

\textbf{Applications} \, One may think of (1) as a version of the fundamental Jacobi field estimates along geodesics used in positive Ricci curvature geometry to derive classical results like the Bonnet-Myers theorem. Similarly, (1)  implies  largeness constraints, now for $scal>0$-geometries. We give some basic examples. Recall, that a compact $n$-manifold  is called \emph{enlargeable} if for any $\ve >0$ there is an orientable covering admitting an $\ve$-contracting map onto the sphere $S^{n}$, i.e. a $C^1$-map with  $\|Df(v)\| \leq \ve \cdot \|v\|$ for all tangent vectors $v$, constant at infinity and of non-zero degree. The Theorem implies extensions of the following result [GL1],Th.12.1, [GL2],Th.A,  [LM], Th.5.5 from the classical to the general case, without changes.\\

\textbf{Corollary 1} \textbf{(Geometry of Enlargeable manifolds)} {\itshape An enlargeable manifold cannot carry $scal > 0$-metrics.
If it carries a $scal \ge 0$-metric, then is covered by a flat torus.}\\

Since the torus $T^n$ is enlargeable, the connected sums $T^n \cs N^n$ with any compact manifold $N^n$  are also enlargeable. Hence, we get the following observation one may paraphrase as follows: there is no general mechanism to locally deform a manifold \emph{increasing} its scalar curvature, even under topological changes. By contrast, we can locally \emph{decrease} the scalar curvature [L1].\\

\textbf{Corollary 2} \textbf{(Non-Existence of $S > 0$-Islands)}   {\itshape There exists \textbf{no} complete Riemannian manifold $(M^n, g)$ such that:
 \begin{itemize}
   \item $scal(g) >0$ on some non-empty open set $U \subset M^n$, with compact closure.
   \item $(M^n \setminus U, g)$ is isometric to $(\R^n \setminus  B_1(0), g_{Eucl})$.
 \end{itemize}}

It is well-known, from [L3], that this result is equivalent to the Riemannian version of the positive mass theorem.\\

\textbf{Corollary 3} \textbf{(Positive Mass Theorem)} {\itshape Let $(M^n,g)$ be asymptotically flat of order $\tau  > \frac{n-2}{2}$. That is, there is a decomposition $M = M_0\cup M_\infty$,
with $M_0$  compact, and $M_\infty$, the end of $M$, diffeomorphic to $\R^n \setminus B_R(0)$, so that this diffeomorphism defines
 coordinates $\{x^i\}$ on $M_\infty$, with {\footnotesize $g_{ij} = \delta_{ij} + O (|x|^{- \tau}), \; \frac{\p g_{ij}}{\p x_k} = O (|x|^{- \tau - 1}),\; \frac{\p^2 g_{ij}}{\p x_k \p x_l} = O(|x|^{- \tau - 2}).$} If $scal(g)\ge 0$, then the \textbf{total mass} {\footnotesize\begin{equation}\label{adm}
E (M,g) := \frac{1}{Vol (S^{n-1})} \cdot \lim_{R \ra
\infty} \int_{\p B_R} \sum_{i,j} \left(\frac{\p g_{ij}}{\p x_i} - \frac{\p g_{ij}}{\p x_j} \right) \cdot \nu_j \: dV_{n-1},
\end{equation}}
where $\nu = (\nu_1\ldots \nu_n)$ is the outer normal vector to $\p B_R$, is \textbf{non-negative}. Moreover, $E (M,g)=0$ if and only if $(M,g)$ is isometric to $(\R^n,g_{Eucl})$.}\\

In a more topological direction, we recall a conjecture related to the Borel conjecture saying that a manifold representing a non-trivial homology class in an \emph{aspherical} manifold  \emph{cannot} carry any $S>0$-metric. This  has been settled in the classical cases cf.[GL1], Th.13.8. As before, these methods immediately extend to the general case. \\

\textbf{ Corollary 4 (Homology of $Sec \le 0$-Manifolds)}{\itshape\, Let $M^n$ be a compact manifold  which admits a $sec \le 0$-metric and $N^k \subset M^n$ be a submanifold with $[N] \neq 0$ in $H_k(M,\Q)$. Then $N^k$ does not admit $scal >0$-metrics.}\\

\textbf{Remark}\, 1.  The present paper gives a more topological argument for the Riemannian positive mass theorem than the one given in  [L4]. Moreover, these methods allow us to simplify our proof of the space-time positive mass theorem [L5] without prior reduction to the Riemannian case or, stated differently, we can cover the Riemannian and the space-time case at once. \\
2. In recent preprints, [CS] and [SY], the authors announce rather different and more classically analytic arguments to address similar results. At any rate, we agree on the statement in [SY] that for such basic results it is valuable to have independent approaches.\\
3. The present paper is based on the potential theory on singular minimal hypersurfaces and their hyperbolic unfoldings [L2]. We use splitting techniques from [L1] and we assume some basic understanding of these methods and their terminology. \qed

As already noted,  the corollaries are well-known consequences of our main theorem discussed, for instance, in [GL1] and  [L3] and we refer to these references for further details. The rest of this paper is devoted to the proof of the Theorem. \\

\setcounter{section}{2}
\renewcommand{\thesubsection}{\thesection}
\subsection{The General Setup} \label{se}

We start with a  compact Riemannian manifold $(M^n,g)=(M^n(\sigma),g(\sigma))$ with boundary  $\partial M^n$, so that $S(g) \ge \sigma$, for some $\sigma>0$  and  some contracting map $f:(M^n, \partial M^n)\rightarrow(S^{n},\{p\})$ of non-zero degree mapping a neighborhood of $\p M^n$ to $p \in S^{n}$. We show inductively that this implies the existence of the following structures  in all dimensions $i$, with $3 \le i \le n$. \\

We get three series of spaces $M^i_{\circledcirc}, M^{i}$ and $P^{i}$ of dimension $i$ so that $M^i_{\circledcirc} \subset M^{i}$, $\partial M^{i} \subset \p M^i_{\circledcirc}$, $M^{n-1}\supset ...  M^{i}  \supset ... M^{3}$ and $M_{\circledcirc}^{n-1}\supset ... M_{\circledcirc}^{i}  \supset ... M_{\circledcirc}^{3}$. They are defined in a way that resembles exact sequences in homological algebra.
 \[M^n_{\circledcirc}  \rhd M^n \rhd P^{n-1} \rhd M^{n-1}_{\circledcirc} \rhd M^{n-1} \rhd \dots P^{i} \rhd M^{i}_{\circledcirc} \rhd M^{i} \rhd \dots P^{3} \rhd M^{3}_{\circledcirc} \]
 where, in our case, $A\rhd B \rhd C$ means that $C$ is determined from $A$ and $B$.\\

In dimension $n$ we choose $M^n_{\circledcirc} = M^n$, $f_n=f$ and $\varphi_n =1$. For $i \le n-1$, the $M_{\circledcirc}^i$ are geometrically regularized subsets of, in general, singular Plateau solutions $P^i$ defined in $M^{i+1}$. In turn, the $M^i$ are smooth extensions of the $M_{\circledcirc}^i$ homologous to the Plateau solution $P^i$. They give mapping degrees for maps onto the sphere $S^i$ a proper sense but there is no control for the curvature on $M^i \setminus M_{\circledcirc}^i$ . In turn, to define the $M^i_{\circledcirc}$ we apply the splitting theorem, [L1], Th.1. \\

\textbf{Definition of $M^{i-1}_{\circledcirc}, M^{i-1}$ and $P^{i-1}$ }: Let $D^{i-1}=\{(x_{1},..., x_{i},0)\in S^{i}: x_{i}\geq0\}$ the lower hemisphere $D^{i-1}\subset S^{i}$ of a geodesic subsphere for $S^{i} \subset \R^{i+1}$ with boundary $\partial D^{i-1}=S^{i-2}$, where
the boundary $\partial D^{i-1}$ can be assumed to be very close to $p_{i}$, but $p_{i} \not \in D^{i-1}$ and  $f_{i}$ is transversal to $D^{i-1}$ and $\partial D^{i-1}$. We get an area minimizing Plateau solution $P^{i-1}$  homologous to $f_{i}^{-1}(D^{i-1})$ in $M^{i}$ for the metric $\varphi_{i} \cdot  g|_{M_{\circledcirc}^{i}}$  and with the boundary $f_{i}^{-1}(S^{i-2})$ which is disjoint to $\p M_{\circledcirc}^{i} \setminus \p M^{i}$. Due to standard boundary regularity results we know that near the smooth boundary $f_{i}^{-1}(S^{i-2})$ the solution $P^{i-1}$  is smooth. Then, and this is where we use the  conformal splitting theorem, [L1], Th.1, we can modify $P^{i-1}$
to define $M_{\circledcirc}^{i-1}$ and, from this, the extension $M^{i-1}$ homologous to $P^{i-1}$. These details are explained in the following chapter.\\

\textbf{Definition of $f_{i-1}$}:  We consider the contracting spherical crushing map $c_i: S^{i}\rightarrow D^{i-1}$ given by
$c_i((x_{1},..., x_{i+1}))=(x_{1},..., x_{i-1},\sqrt{x_{i}^{2}+x_{i+1}^{2}},0), \mm{ for } (x_{1},...,x_{i+1}) \in S^{i} \subset \R^{i+1}.$ We compose $c_i$ with a  map $d_i: D^{i-1}\rightarrow S^{i-1}$, which maps the boundary $\partial D_{i-1}$ and $c(p_{i})$ to a point $p_{i-1}\in S^{i-1}$ to define a  $4$-contracting map $h_{i}:=d_i \circ c_i: S^{i}\rightarrow S^{i-1}$. From the construction of the map $f_{i-1}$ it follows that these new boundary components are all mapped onto the point $p_{i-1}\in S^{i-1}$.  We define $f_{i-1}:=h_{i} \circ f_{i}|_{M^{i-1}}$. It follows that $f_{i-1}$ is $a_{i-1}$-contracting, where $a_i$ only depends on $i$.  We have $deg(f_{i-1})\neq 0$ since $deg (h_{i}|_{D_{i-1}})\neq 0$, $deg (f_{i})\neq 0$ and $N$ being homologous to $f_{i}^{-1}(D^{i-1})$.\\

\textbf{Definition of $\varphi_{i-1}$}:
There are smooth functions $\varphi_{i-1} >0$ on the cores $M_{\circledcirc}^{i-1}$  with $\varphi_{i-1}\cdot  S(\varphi_{i-1}\cdot g_{M^{i-1}}) \ge\sigma$ and the $M_{\circledcirc}^{i}$ are regular subsets of the generally singular Plateau solution $P^{i-1}$ with $\p P^{i-1}=f_{i}^{-1}(S^{i-2})$ relative the metric $\varphi_{i+1}\cdot g_{M^{i+1}}$. The definition of the $\varphi_{i-1}$ is related to that of $M_i$ and it will be postponed to chapter 4. The weaker  condition  $\varphi_{i}\cdot  S(\varphi_{i}\cdot g_{M^{i}}) \ge\sigma$, compared to  $S(g_{M^{n}}) \ge\sigma$, has the advantage to survive the dimensional descent we use here.\\

\textbf{Final Conclusion and Summary}:  We run this process until we reach dimension $3$. A simple adjustment allows us to continue to dimension $1$ where we argue using ordinary differential inequalities to find an upper bound for $\sigma$ depending only on $n$. The overall setup deviates from the classical one in  [GL1],Ch.12 in several ways. We incorporate singular Plateau solutions $P^i$ and regularize $P^i$, in a geometric step, to $M^i_{\circledcirc}$ and, in a topological one, to $M^i$. The step from $P^i$ to $M^i_{\circledcirc}$ is not compatible with the symmetrization procedure in  [GL1],Ch.12 and we replace it, even in low dimensions,  for a direct conformal deformation.\\

Henceforth, we denote the collection of the $i$ dimensional data $M^i_{\circledcirc}, M^{i}$ and $P^{i}$ with the maps $f_i$ and $\varphi_i$ by $\mathbf{C}_i$ and the transition from  $\mathbf{C}_i$ to $\mathbf{C}_{i-1}$ by $\rhd$.

\setcounter{section}{3}
\renewcommand{\thesubsection}{\thesection}
\subsection{Conformal Splittings and $\mathbf{C}_i \rhd \mathbf{C}_{i-1}$} \label{se}

We start with the basic step $\mathbf{C}_n \rhd \mathbf{C}_{n-1}$.  For the given compact Riemannian manifold $(M^n,g)$ with $S(g) \ge \sigma$  and the contracting map $f:(M^n, \partial M^n)\rightarrow(S^{n},\{p\})$ we choose a lower hemisphere $D^{n-1}\subset S^{n}$ with boundary $\partial D^{n-1}=S^{n-2}$, so that $\partial D^{n-1}$ is close to $p$, but $p \not \in D^{n-1}$ and  we may assume that $f$ is transversal to $D^{n-1}$ and $\partial D^{n-1}$.\\

We get an area minimizing Plateau solution $P^{n-1}$  homologous to $f^{-1}(D^{n-1})$ in $M^{n}$  with the boundary $f^{-1}(S^{n-2})$.
It may happen that $P^{n-1}$ hits $\p M$, but this does not cause problems since $f_{n-1}$ maps also a neighborhood of $P^{n-1} \cap \p M$   the new basepoint $p_{n-1}$. Thus, in our present argument we can assume that  $P^{n-1} \cap \p M \subset \p P^{n-1}$. (In these places the regularity of $P$ is not needed.) Since $P^{n-1}$ is area minimizing we have $A''(f) := Area''(f \cdot \nu) \ge 0$ for any smooth infinitesimal variation supported away from the singular set $\Sigma_P \subset P$ and the boundary $f \cdot \nu$ of $P \setminus (\Sigma_P \cup \p P)$, where $\nu$ is a unit normal vector field, (we may
assume $M$ is orientable). A direct computation shows:

{\small \[A'' (f) = \int_P | \nabla_P f |^2 - f^2 (| A |^2 + Ric_M (\nu, \nu)) d A \ge 0 \; \Longleftrightarrow \]
\[\int_P | \nabla f |^2 + \frac{n-2}{4 (n-1)} (scal_P-scal_M) f^2 d A \ge \int_P \frac{n}{2 (n-1)} | \nabla f |^2 + \frac{n- 2}{4 (n-1)} f^2  | A |^2   d A.\]}
Thus for a Hardy \si-transform $\bp$ we get
{\small \[\int_P | \nabla f |^2 + \frac{n-2}{4 (n-1)} (scal_P-scal_M) f^2 d A \ge \lambda_{\bp} \cdot \int_P \bp^2\cdot f^2 dA,\]}
for some (largest) constant $\lambda_{\bp}>0$. In particular  we get on $P\setminus \Sigma_P$ the \emph{ground state} for the \si-adapted Schr\"odinger operator $L(P):=-\frac{4(n-1)}{n-2}\triangle+scal(g|_P) -scal_M$, that is, a smooth function $\psi>0$ with $\psi \equiv 0$ on $\p P$  and with
{\small
\begin{equation}\label{sca}
-\frac{4(n-1)}{n-2}\triangle\psi+scal(g|_P) \psi -scal_M \psi=  \lambda_{\bp}\langle A\rangle^{2} \psi
\end{equation}
}

This observation and growth estimates for such  functions $\psi$ from [L2], Th.6 are the basis for the following result which is a version of [L1], Th.1 in the case $i=n-1$. The cases $3 \le i <n-1$, where we first need to define appropriate Plateau solutions will be proved inductively.

\begin{proposition}\emph{\textbf{(Conformal Splitting in Dimension $i$)}}\label{cs} For any $\ve >0$, there are smooth functions $\lambda_\ve, \phi_{P^i}>0$ on $P^i\setminus \Sigma_{P^i}$  with $\phi_{P^i} \equiv 0$ on $\p {P^i}$, with $c_{P^i} \le \lambda_{\ve} \le C_{P^i}$, for some constants $c_{P^i}, C_{P^i}>0$ and a neighborhood $U_{\ve}$ of $\Sigma_{P^i}$, within an $\ve$-neighborhood of $\Sigma_{P^i}$ so that
{\small\begin{equation}\label{la}
L_{\ve}:=-\frac{4(i-1)}{i-2}\triangle +scal(g|_{P^i})-scal_M- \lambda_{\ve} \cdot \langle A\rangle^{2}\mm{ is \si-adapted and } L_{\ve} \phi_{P^i}=0
\end{equation}}

and $M^i_{\circledcirc}= M^i_{\circledcirc}(\ve)=(P^i \setminus U_{\ve},\phi^{4/i-2}_{P^i} \cdot g_P)$ is a $scal>0$-manifold with (locally) area minimizing boundary $\p M^i_{\circledcirc}$. Moreover, there is a smooth manifold $M^i$ with boundary $\p P^i$ with $M^i_{\circledcirc} \subset M^i$  homologous to $P^i$.
\end{proposition}

\begin{remark} In the following argument $\ve>0$ is a parameter we gradually shrink to a tiny value.  For $i \le 7$,  we choose  $U_{\ve} \v$ and $M^i=M^i_{\circledcirc}=P^i$. For $i \ge 8$ the extension from $M^i_{\circledcirc}$  to $M^i$  is the relative version of the general observation that for any closed manifold $X^n$ one has $H_{n-1}(X,\Z) \cong H^{1}(X,\Z)  \cong [X,S^1]$. Namely, for  any $\alpha \in H_{n-1}(X,\Z)$ and any smooth map $f: X \ra S^1$ that represents the homotopy class in $[X,S^1]$ associated to $\alpha$, the smooth preimage $f^{-1}(c)$ of some generic $c \in S^1$ represents $\alpha$.\qed
\end{remark}

For the induction step $\mathbf{C}_i \rhd \mathbf{C}_{i-1}$, $i  <n$, we start from $M^i$. We observe that \ref{cs} does not control for the curvatures on $M^i \setminus M^i_{\circledcirc}$. However, the area minimizing Plateau solution $P^{i-1}$  homologous to $f_{i-1}^{-1}(D^{i-1})$ in $M^i$  with the boundary $f_{i-1}^{-1}(S^{i-2})$ may intersect this subset. We use the minimality of
$\p M^i_{\circledcirc}$ to rebuild $P^{i-1}$ to bypass $M^i \setminus M^i_{\circledcirc}$ and to define the next two spaces $M^{i-1}_{\circledcirc}$ and  $M^{i-1}$.\\

To this end, we may assume, after selecting one of its components, that $\p M^i_{\circledcirc}$ is connected.  For small $\ve >0$, we may assume that $\p P^{i-1} \cap \p M^i_{\circledcirc} \v$. If  $P^{i-1} \cap \p M^i_{\circledcirc} \n$ the intersection is transversal in all regular points due to the strict maximum principle and there are two cases where $\p M^i_{\circledcirc}$ plays the role of an either \emph{repelling} or \emph{attracting} end:\\

\textbf{Repelling Ends}\,  $\p P^{i-1} \cap \p M^i_{\circledcirc}= \p U$ for some open subset $U \subset \p M^i_{\circledcirc}$. Then we can replace $(P^{i-1} \setminus (M^i \setminus M^i_{\circledcirc}))\cup U$ for some other Plateau solution $Q^{i-1}\cap \p M^i_{\circledcirc} \v$ again using the strict maximum principle.\\

\textbf{Attracting Ends}\,  Otherwise,  we may assume that $\p P^{i-1} \cap \p M^i_{\circledcirc}$ locally separates $\p M^i_{\circledcirc}$ but  $\p M^i_{\circledcirc} \setminus (P^{i-1} \cap \p M^i_{\circledcirc})$ remains connected.
  Then, for small $\ve >0$ and any $m \in \Z^{\ge 1}$, we get, using the strict maximum principle another time:
 \begin{itemize}
   \item an at least locally  area minimizing $Q_m^{i-1}$ with two (not necessarily connected) disjoint boundary components $\p Q_m^{i-1}[A]:=\p P^{i-1}$ and the \emph{far end}
$\p Q_m^{i-1}[B]:=P^{i-1}\cap \p M^i_{\circledcirc}$.
   \item $[Q_m^{i-1}]$ is homologous to $[P^{i-1} \cap M^i_{\circledcirc}] + m \cdot [(\p M^i_{\circledcirc} \setminus (P^{i-1}\cap \p M^i_{\circledcirc}))^\sim]$ in $H_{n-1}(M^i_{\circledcirc} , \p M^i_{\circledcirc};\Z)$, where  $\sim$ means the non-embedded one-sided closure of $\p M^i_{\circledcirc} \setminus (P^{i-1}\cap \p M^i_{\circledcirc})$.
 \end{itemize}

One may  compare this case with that of a geodesic ray on a surface approaching a closed geodesic in a gradually narrowing spiral.\\

We observe, for $\ve >0$ small and $m$ large enough, starting and counting backwards from $\p Q_m^{i-1}[B]$, the hypersurface $Q_m^{i-1}$  can be chopped into  $m$ pieces, we choose $m$ even, $V_{m,k} \subset Q_m^{i-1}$, $k=1,..m$, so that viewed as currents
\begin{equation}\label{asyy}
V_{m,k} \ra \p M^i_{\circledcirc} \setminus (P^{i-1}\cap \p M^i_{\circledcirc}) \mm{ for } k \le m/2, m \ra \infty, \mm{ in flat norm.}
\end{equation}
We apply \ref{cs} to the possibly singular $\p M^i_{\circledcirc}$ in its $scal>0$-ambient space $M^i_{\circledcirc}$ and get, as an intermediate step, another conformal $scal >0$-geometry written $N_{\circledcirc}(\p M^i_{\circledcirc})^{i-1}$ with area minimizing boundary $\p N_{\circledcirc}^{i-2}$ to some intermediate smooth manifold $N^{i-1}$. \\

Next we also apply the splitting result  to $Q_m^{i-1}$  and notice from the use of the naturality of the underlying Martin theory for the family of Schr\"odinger operators $L(\cdot)$ form  [L2], Th.3 to periodically define the same $scal >0$-geometry and extensions for the $V_{m,k}$. We call the $scal>0$-geometry $W_{m,_{\circledcirc}}^{i-1}$ with minimal border $\p W_{m,{\circledcirc}}^{i-1}$ and the smooth extension geometry $W_m^{i-1}$, we can  define, from (\ref{asyy}), to become  asymptotically periodical for large $m$ and $k \le m/2$. The limit for $m \ra \infty$ is written $W_{{\circledcirc}}^{i-1}$  and asymptotically periodical  extension pieces we associate to the $V_{m,k}$ are labelled  $WV_{m,k}$.\\

Next, for sufficiently large $m \gg 1$ we may assume that we have an area minimizing hypersurface $X^{i-2} \subset WV_{m,k-1} \cup ... \cup WV_{m,k-l} \subset Q_m^{i-1}$, for $k = 1/2 \cdot m , l  = 1/4 \cdot m$, separating  $WV_{m,k}$ from $WV_{m,k-l}$.   This can easily be accomplished, e.g. applying a slight scaling of the pieces of $W_m^{i-1}$ towards infinity, which keeps the scalar curvature condition (\ref{n})  below up to a factor $\in (1/2,2)$. \\

We write $W_{m,X}^{i-1}$ for the component that contains  $WV_{m,k}$.  Again we are in the situation that $X^{i-2}$ may intersect $\p W_{m,_{\circledcirc}}^{i-1}$ and reach a subset without controlled curvatures.
However, in this case we are in an appended end piece of $Q_m^{i-1}$. Hence, for $\ve >0$ small enough and large $m$, the boundary $\p W_{m,X}^{i-1}$ is a repelling end. Thus we get other area minimizing boundaries $Y^{i-2}=Y^{i-2}(m,\ve) \subset W_{m,X}^{i-1}$ with $Y^{i-2}(m,\ve) \ra X^{i-2}$ for $\ve \ra 0$ and $Y^{i-2}(m,\ve) \ra W_{{\circledcirc}}^{i-1}$  for $m \ra \infty$, in flat norm.  We first choose a large $m$ and then shrink $\ve$ to define the new $M^{i-1}_{\circledcirc}$ as the component of $W_{{\circledcirc}}^{i-1}$ that does \emph{not} contain $X^{i-2}$. This reproduces \ref{cs} in dimension $i-1$ and from this we readily get the other components of $\mathbf{C}_{i-1}$.\\

\setcounter{section}{4}
\renewcommand{\thesubsection}{\thesection}
\subsection{Recursive Conformal Deformations} \label{p2}

Now we discuss the propagation of  scalar curvature estimates while we run through the steps $\mathbf{C}_i \rhd \mathbf{C}_{i-1}$. As before we start from  be a compact manifold $M^n$ with boundary  $\partial M^n$ and a  contracting $C^1$-map $f:(M^n, \partial M^n)\rightarrow(S^{n},\{p\})$ of non-zero
degree, which maps a neighborhood of $\p M^n$ to a point $p \in S^{n}$. The compactness gives some $\kappa \in \R$ so that $scal(g_M) \ge \kappa$ and our goal is to show that there is upper bound on the possible values of $\kappa$ for any such given $M$, $f$ depending only on the dimension $n$. \\

 We recall the scalar curvature transformation law for an $k$-dimensional Riemannian manifold $(N^k, g_N)$, $k \ge 3$, under conformal deformations $u^{4/k-2} \cdot g_N$ for smooth $u>0$:
{\small \begin{equation} \label{4} scal(u^{4/k-2} \cdot g_N) \cdot u^{\frac{k+2}{k-2}} =  -\Delta_N u +\frac{k-2}{4 (k-1)} \cdot scal_N \cdot u
\end{equation}}
From this equality and $scal(g_M) \ge \kappa$ we get for a positive solutions $\psi_{M_{\circledcirc}^{n-1}}$ of (\ref{sca}):
{\small \begin{equation}\label{n}
 scal(\psi_{M_{\circledcirc}^{n-1}}^{4/n-3} \cdot  g_{M_{\circledcirc}^{n-1}}) \cdot \psi_{M_{\circledcirc}^{n-1}}^{4/n-3} \ge \kappa.
 \end{equation}}
Now we argue inductively in two steps. We choose the area minimizer of the last chapter and find for some positive functions $\lambda_{[i]}>0$ from \ref{cs}, for some small $\ve >0$, so that
{\small\[(M_{\circledcirc}^{n-2}, \psi_{M_{\circledcirc}^{n-1}}^{4/n-3} \cdot  g_{M_{\circledcirc}^{n-2}})\subset (M_{\circledcirc}^{n-1},\psi_{M_{\circledcirc}^{n-1}}^{4/n-3}\cdot  g_{M_{\circledcirc}^{n-1}})\]
\[scal\left(\psi_{M_{\circledcirc}^{n-2}}^{4/n-4} \cdot (\psi_{M_{\circledcirc}^{n-1}}^{4/n-3} \cdot g_{M_{\circledcirc}^{n-2}})\right) \cdot \psi_{M_{\circledcirc}^{n-2}}^{4/n-4}\ge
 \lambda_{[n-2]} \cdot \langle A\rangle_{M_{\circledcirc}^{n-2}}^{2}   + scal(\psi_{M_{\circledcirc}^{n-1}}^{4/n-3}\cdot g_{M_{\circledcirc}^{n-1}}) \]}

We multiply by $\psi_{M_{\circledcirc}^{n-1}}^{4/n-3}$ to apply the inductive assumption

{\small\begin{equation}\label{n-1}
scal\left(\psi_{M_{\circledcirc}^{n-2}}^{4/n-4} \cdot (\psi_{M_{\circledcirc}^{n-1}}^{4/n-3} \cdot g_{M_{\circledcirc}^{n-2}})\right) \cdot \psi_{M_{\circledcirc}^{n-2}}^{4/n-4} \cdot \psi_{M_{\circledcirc}^{n-1}}^{4/n-3} \ge
\end{equation}
\[\lambda_{[n-2]} \cdot \langle A\rangle_{M_{\circledcirc}^{n-2}}^{2} \cdot \psi_{M_{\circledcirc}^{n-1}}^{4/n-3}  + scal(\psi_{M_{\circledcirc}^{n-1}}^{4/n-3}\cdot g_{M_{\circledcirc}^{n-1}})\cdot \psi_{M_{\circledcirc}^{n-1}}^{4/n-3} \ge  \lambda_{[n-2]} \cdot \langle A\rangle_{M_{\circledcirc}^{n-2}}^{2} \cdot \psi_{M_{\circledcirc}^{n-1}}^{4/n-3} +\kappa \ge \kappa.\]}
Next we choose
{\small\[(M_{\circledcirc}^{n-3}, \psi_{M_{\circledcirc}^{n-2}}^{4/n-4} \cdot (\psi_{M_{\circledcirc}^{n-1}}^{4/n-3} \cdot g_{M_{\circledcirc}^{n-3}}))\subset (M_{\circledcirc}^{n-2},\psi_{M_{\circledcirc}^{n-2}}^{4/n-4} \cdot (\psi_{M_{\circledcirc}^{n-1}}^{4/n-3} \cdot g_{M_{\circledcirc}^{n-2}}))\]}
and repeat the previous arguments iteratively until we reach
{\small \begin{equation}\label{3}
scal\left(\psi_{M_{\circledcirc}^{3}}^{4} \cdot .... \cdot \psi_{M_{\circledcirc}^{n-2}}^{4/n-4} \cdot \psi_{M_{\circledcirc}^{n-1}}^{4/n-3} \cdot g_{M_{\circledcirc}^{3}}\right) \cdot \psi_{M_{\circledcirc}^{3}}^{4} \cdot ... \cdot \psi_{M_{\circledcirc}^{n-2}}^{4/n-4} \cdot \psi_{M_{\circledcirc}^{n-1}}^{4/n-3} \ge
\end{equation}
\[\lambda_{[3]} \cdot \langle A\rangle_{M_{\circledcirc}^{3}}^{2} \cdot \psi_{M_{\circledcirc}^{4}}^{4/2}  ... \cdot \psi_{M_{\circledcirc}^{n-2}}^{4/n-4} \cdot \psi_{M_{\circledcirc}^{n-1}}^{4/n-3} + scal(\psi_{M_{\circledcirc}^{4}}^{4/2}  ... \cdot \psi_{M_{\circledcirc}^{n-2}}^{4/n-4} \cdot \psi_{M_{\circledcirc}^{n-1}}^{4/n-3} \cdot g_{M_{\circledcirc}^4})\cdot \psi_{M_{\circledcirc}^{4}}^{4/2}  ... \cdot \psi_{M_{\circledcirc}^{n-2}}^{4/n-4} \cdot \psi_{M_{\circledcirc}^{n-1}}^{4/n-3}\]
  \[\ge  \lambda_{[3]} \cdot \langle A\rangle_{M_{\circledcirc}^{3}}^{2} \cdot \psi_{M_{\circledcirc}^{4}}^{4/2}  ... \cdot \psi_{M_{\circledcirc}^{n-2}}^{4/n-4} \cdot \psi_{M_{\circledcirc}^{n-1}}^{4/n-3} +\kappa \ge \kappa.\]}

To extend this dimensional descent to the two and finally the one dimensional case we observe that the argument for [L1], Th.1 and the estimates [L2], Th.6 (and specifically Prop.4.6) show that we can apply the previous arguments not only to $M$ but also to the Riemannian product $M \times S^1 \times S^1$ for conformal deformations of $M \times S^1 \times S^1$ depending only on the $M$. This makes the various deformation and construction steps we have seen before invariant under the $S^1 \times S^1$-action and we finally reach a $3$-manifold with connected components of the form $I \times S^1 \times S^1$, where $I \subset \R$ is an interval.

\begin{remark} We note in passing that the bending effect towards singularities, needed to find an area minimizing border, would not be strong enough to make the symmetrization argument of [GL], Ch.12 work in high dimension.\qed
\end{remark}

For such a \emph{flat} three manifold $I \times S^1 \times S^1$ we have from (\ref{la}), (\ref{4}) and (\ref{3})
{\small \begin{equation}\label{33}
-8\frac{\triangle \psi}{\psi}= scal(\varphi \cdot(g|_{I}+g_{S^{1}\times S^{1}})) \cdot \varphi\ge \kappa \mm{  with } \psi^{4}:=\varphi.
\end{equation}}
For at least one component  $I \times S^1 \times S^1$ we have  $f_{1}:I\rightarrow S^{1}$ is of nonzero degree with $f_{1}(I)=S^{1}$. After reparametrization of $I$ by arc-length we have that $f_1$ is $4^n$-contracting and this can be reduced to the condition $\frac{\overline{\psi}^{''}}{\overline{\psi}}\ge \kappa/8$ on the interval $[0,4^{-n}]$, where $\overline{\psi}(x)=\int_{\{x\} \times S^1 \times S^1}\psi \, dV$, since $\int_{\{x\} \times S^1 \times S^1} \Delta_{\{x\} \times S^1 \times S^1}\psi|_{\{x\} \times S^1 \times S^1} \, dV=0$.
 An elementary computation, very similar to that in [GL], Ch.12, then gives the desired upper bound on $\kappa$.\qed

\end{document}